\documentclass[12pt]{amsart}

\usepackage{amsmath,xspace,amssymb,mathrsfs}
\usepackage{color}

\input xy
\xyoption{all}
\xyoption{2cell}
\UseAllTwocells
\CompileMatrices

\newcommand{\Spec}{\operatorname{Spec}}
\renewcommand{\phi}{\varphi}

\newcommand{\Ker}{\operatorname{Ker}}
\newcommand{\Ima}{\operatorname{Im}}

\newcommand{\Coker}{\operatorname{Coker}}

\newtheorem{proposition}{Proposition}[section]
\newtheorem{lemma}[proposition]{Lemma}

\newtheorem{theorem}[proposition]{Theorem}

\theoremstyle{definition}
\newtheorem{definition}[proposition]{Definition}

\newtheorem{remark}[proposition]{Remark}

\begin{document}

\title{On finite type epimorphisms of rings}

\author{Abolfazl Tarizadeh}
\address{Department of Mathematics, Faculty of Basic Sciences, University of Maragheh \\
P. O. Box 55136-553, Maragheh, Iran.
 }
\email{ebulfez1978@gmail.com}

\date{}
\footnotetext{ 2010 Mathematics Subject Classification: 13B10, 13B30, 13A99, 13C11. \\ Key words and phrases: finite type epimorphism; module of differentials; universally injective.}

\begin{abstract} In this note, finite type epimorphisms of rings are  characterized.
\end{abstract}

\maketitle

\section{Introduction}

In this paper we investigate an important class of epimorphisms of rings, namely finite type epimorphisms. Needless to say that finite type ring maps  have very closed connections with the geometric structures and essentially originate from algebraic geometry. Meanwhile, finite type epimorphisms of commutative rings have very interesting properties and they are important from various aspects. For example, ``finite type monomorphisms of schemes'' can be considered as the geometric interpretation of finite type epimorphisms of rings. \\

In the realm of epimorphisms of commutative rings there are some highly non-trivial results in the literature. For instance, ``every finite type epimorphism of rings which is also injective and flat then it is of finite presentation'', see \cite[Theorem 1.1]{Cox-Rush}. A special case of this result was announced in \cite[Corollary 3.4.7]{Raynaud-Gruson}. Another important result is due to Lazard, see Theorem \ref{Theorem I}. In the present article, partially motivated by the above results, we have obtained two new and non-trivial results. In fact, Theorems \ref{Theorem II} and \ref{geo th} are the main results of this note. These results can be considered as the analogous of Theorem \ref{Theorem I} in the finite type case whose hypotheses have been as much as possible weakened. In this  article, all of the rings are commutative. \\

\section{preliminaries}

Here we recall some material which are needed in the next section. \\

By an epimorphism $\phi:R\rightarrow S$ we mean it is an epimorphism in the category of (commutative) rings. Surjective ring maps are special cases of epimorphisms. As a specific example, the canonical ring map $\mathbb{Z}\rightarrow\mathbb{Q}$ is an epimorphism of rings while it is not surjective.  \\

\begin{remark}\label{Remark II} It is easy to see that a ring map $\phi:R\rightarrow S$ is an epimorphism of rings if and only if $s\otimes1=1\otimes s$ for all $s\in S$. It is also equivalent to the condition that $S\otimes_{R}\Coker\phi=0$. It follows that every faithfully flat epimorphism of rings is an isomorphism. In particular, every non-zero epimorphism of rings with source a field is an isomorphism. We refer to \cite{Samuel}, specially \cite{Lazard}, \cite{Olivier} and \cite{Roby} for a comprehensive discussion of epimorphisms of commutative rings. Also see \cite{Call}. \\
\end{remark}

\begin{theorem}\label{Theorem I} A ring map $\phi: R\rightarrow S$ is an epimorphism of rings if and only if the following conditions hold.\\
$\mathbf{(i)}$ The induced map $\phi^{\ast}:\Spec(S)\rightarrow\Spec(R)$ is injective.\\
$\mathbf{(ii)}$ For each prime ideal $\mathfrak{q}$ of $S$ the induced map $\kappa(\mathfrak{p})\rightarrow\kappa(\mathfrak{q})$ is an isomorphism where $\mathfrak{p}=\phi^{\ast}(\mathfrak{q})$.\\
$\mathbf{(iii)}$ The kernel of the canonical ring map $S\otimes_{R}S\rightarrow S$ is a finitely generated ideal.\\
$\mathbf{(iv)}$ The module of K\"{a}hler differentials  $\Omega_{S/R}$ is zero.\\
\end{theorem}

{\bf Proof.} See \cite[Proposition 1.5]{Lazard}. $\Box$ \\

\begin{remark}\label{Remark I} It is well-known that if $\phi:R\rightarrow A$ is a ring map and $\mathfrak{p}$ a prime ideal of $R$ then $\mathfrak{p}\in\Ima\phi^{\ast}$ if and only if $A\otimes_{R}\kappa(\mathfrak{p})\neq0$ where $\kappa(\mathfrak{p})$ is the residue field of $R$ at $\mathfrak{p}$. \\
\end{remark}

\begin{definition} The map $\phi^{\ast}:\Spec(S)\rightarrow\Spec(R)$ induced by a ring map $\phi:R\rightarrow S$ is said to be \emph{universally injective} if for any ring map $R\rightarrow R'$ then the induced map $\psi^{\ast}:\Spec(R'\otimes_{R}S)\rightarrow\Spec(R)$ is injective where $\psi: R'\rightarrow R'\otimes_{R}S$ is the base change map.\\
\end{definition}

\section{Main results}

Let $\phi:R\rightarrow S$ be a ring map. Consider the canonical ring map $\pi: S\otimes_{R}S\rightarrow S$ which maps each pure tensor $s\otimes s'$ of $S\otimes_{R}S$ into $ss'$. The kernel of $\pi$ is generated by elements of the form $s\otimes1-1\otimes s$. Because if $\sum\limits_{i}s_{i}s'_{i}=0$ then we may write $\sum\limits_{i}s_{i}\otimes s'_{i}=\sum\limits_{i}(1\otimes s'_{i})(s_{i}\otimes1-1\otimes s_{i})$. \\

To prove the main results of this note we need the following lemma. \\

\begin{lemma}\label{lemma 1} If a ring map $\phi:R\rightarrow S$ is of finite type then $\Ker\pi$ is a finitely generated ideal.\\
\end{lemma}

{\bf Proof.} By the hypothesis there are elements $s_{1},...,s_{n}\in S$ such that $S=R[s_{1},...,s_{n}]$. Consider the ideal $J=(s_{i}\otimes1-1\otimes s_{i}: 1\leq i\leq n)$ of $S\otimes_{R}S$. Clearly $J\subseteq\Ker\pi$. To prove the reverse inclusion it suffices to show that  $$s_{1}^{d_{1}}...s_{n}^{d_{n}}\otimes1-1\otimes s_{1}^{d_{1}}...s_{n}^{d_{n}}\in J.$$ We use an induction argument over $n$. If $n=1$ then we have $s^{d}\otimes1-1\otimes s^{d}=(s\otimes1)^{d}-(1\otimes s)^{d}=\big((s\otimes1)^{d-1}+(s\otimes1)^{d-2}(1\otimes s)+...+(1\otimes s)^{d-1}\big)(s\otimes1-1\otimes s)$ which belongs to $J$. Let $n>1$. Then we may write $s_{1}^{d_{1}}...s_{n}^{d_{n}}\otimes1-
1\otimes s_{1}^{d_{1}}...s_{n}^{d_{n}}=
s_{n}^{d_{n}}\otimes1\big(s_{1}^{d_{1}}...s_{n-1}^{d_{n-1}}\otimes1-
1\otimes s_{1}^{d_{1}}...s_{n-1}^{d_{n-1}}\big)+1\otimes s_{1}^{d_{1}}...s_{n-1}^{d_{n-1}}\big(s_{n}^{d_{n}}\otimes1-1\otimes s_{n}^{d_{n}}\big)$ which is, by the induction hypothesis and the induction step, belonging to $J$. $\Box$ \\

\begin{theorem}\label{Theorem II} A finite type ring map $\phi:R\rightarrow S$ is an epimorphism of rings if and only if the induced map $\phi^{\ast}$ is injective and for each prime ideal $\mathfrak{q}$ of $S$
the base change map $\kappa(\mathfrak{p})
\rightarrow S_{\mathfrak{q}}\otimes_{R}\kappa(\mathfrak{p})$
is an isomorphism where $\mathfrak{p}=\phi^{\ast}(\mathfrak{q})$. \\
\end{theorem}

{\bf Proof.} ``$\Rightarrow$'' By Theorem \ref{Theorem I}, $\phi^{\ast}$ is injective. The composition of $\phi$ with the canonical ring map $S\rightarrow S_{\mathfrak{q}}$ gives us the epimorphism $R\rightarrow S_{\mathfrak{q}}$. Thus the map $\kappa(\mathfrak{p})
\rightarrow S_{\mathfrak{q}}\otimes_{R}\kappa(\mathfrak{p})$ is an epimorphism since every epimorphism of rings is stable under the base change. By Remark \ref{Remark I}, $S_{\mathfrak{q}}\otimes_{R}\kappa(\mathfrak{p})\neq0$ and so by Remark \ref{Remark II} the base change map is an isomorphism. To prove the converse implication we shall use Theorem \ref{Theorem I}. We have $S_{\mathfrak{q}}\otimes_{R}\kappa(\mathfrak{p})\simeq S_{\mathfrak{q}}/\mathfrak{p}S_{\mathfrak{q}}$. Thus $S_{\mathfrak{q}}/\mathfrak{p}S_{\mathfrak{q}}$ is a field and so $\mathfrak{p}S_{\mathfrak{q}}=\mathfrak{q}S_{\mathfrak{q}}$. Hence the induced map $\kappa(\mathfrak{p})\rightarrow\kappa(\mathfrak{q})$ is an isomorphism. Moreover, by \cite[Tag 00RV]{Johan}, we have
$$\Omega_{S_{\mathfrak{q}}/R}\otimes_{R}\kappa(\mathfrak{p})
\simeq\Omega_{S_{\mathfrak{q}}
\otimes_{R}\kappa(\mathfrak{p})/\kappa(\mathfrak{p})}=0.$$
It follows that
$(\mathfrak{q}S_{\mathfrak{q}})
\Omega_{S_{\mathfrak{q}}/R}=
\Omega_{S_{\mathfrak{q}}/R}$.
By \cite[Tag 00RZ]{Johan}, $\Omega_{S/R}$ is finitely generated $S-$module.
Therefore $(\Omega_{S/R})_{\mathfrak{q}}\simeq\Omega_{S_{\mathfrak{q}}/R}$ is a finitely generated $S_{\mathfrak{q}}-$module. Then Nakayama implies that $\Omega_{S_{\mathfrak{q}}/R}=0$. It follows that $\Omega_{S/R}=0$
because $(\Omega_{S/R})_{\mathfrak{q}}\simeq\Omega_{S_{\mathfrak{q}}/R}=0$.  Now using Lemma \ref{lemma 1} and Theorem \ref{Theorem I} then we conclude that $\phi$ is an epimorphism.  $\Box$ \\

\begin{theorem}\label{geo th} Let $\phi:R\rightarrow S$ be a finite type ring map. Then the following conditions are equivalent.\\
$\mathbf{(i)}$ $\phi$ is an epimorphism of rings.\\
$\mathbf{(ii)}$  $\phi^{\ast}$ is universally injective and  $\Omega_{S/R}=0$.\\
$\mathbf{(iii)}$ $\phi$ is formally unramified and for any field $K$ and any ring maps $f,g:S\rightarrow K$ if $f\circ\phi=g\circ\phi$ then $f=g$.\\
$\mathbf{(iv)}$  $\phi^{\ast}$ is injective, $\Omega_{S/R}=0$ and for each prime ideal $\mathfrak{q}$ of $S$ the field extension $\kappa(\mathfrak{p})\subseteq\kappa(\mathfrak{q})$ is purely inseparable where $\mathfrak{p}=\phi^{\ast}(\mathfrak{q})$.\\
\end{theorem}

{\bf Proof.} $\mathbf{(i)}\Rightarrow\mathbf{(ii)}:$ Epimorphisms are stable under the base change. Then apply Theorem \ref{Theorem I}.\\
$\mathbf{(ii)}\Rightarrow\mathbf{(iii)}:$ The ring map $\phi$ is formally unramified if and only if $\Omega_{S/R}=0$, see \cite[Tag 00UO]{Johan}. There exist ring maps $f',g':K\otimes_{R}S\rightarrow K$ which map each pure tensor $a\otimes s$ of $K\otimes_{R}S$ into $af(s)$ and $ag(s)$, respectively. Let $\psi:K\rightarrow K\otimes_{R}S$ be the base change of  $f\circ\phi:R\rightarrow K$. Then, by the hypotheses, $\psi^{\ast}$ is injective.
Moreover $K\otimes_{R}S$ is a nontrivial ring. Therefore the prime spectrum of $K\otimes_{R}S$ is a single-point set. This implies that $(f')^{-1}(0)=(g')^{-1}(0)$. But for each $s\in S$, $f'\big(f(s)\otimes1-1\otimes s\big)=0$. Thus $g'\big(f(s)\otimes1-1\otimes s\big)=0$. Therefore $f(s)=g(s)$ for all $s\in S$.\\
$\mathbf{(iii)}\Rightarrow\mathbf{(i)}:$ Let $\mathfrak{p}$ be a prime ideal of $S\otimes_{R}S$. The hypotheses imply that
$\eta\circ i=\eta\circ j$ where $\eta:S\otimes_{R}S\rightarrow\kappa(\mathfrak{p})$ and $i,j:S\rightarrow S\otimes_{R}S$ are the canonical ring maps. Therefore for each $s\in S$, $s\otimes1-1\otimes s$ is a nilpotent element. Let $J$ be the kernel of the canonical ring map $S\otimes_{R}S\rightarrow S$.
By Lemma \ref{lemma 1}, $J$ is a finitely generated ideal.
It follows that $J$ is a nilpotent ideal. But $\Omega_{S/R}\simeq J/J^{2}$. Thus $J=J^{2}$ and so $J=0$. Therefore $s\otimes1=1\otimes s$ for all $s\in S$. Hence $\phi$ is an epimorphism. $\mathbf{(i)}\Rightarrow\mathbf{(iv)}:$ See Theorem \ref{Theorem I}.\\
$\mathbf{(iv)}\Rightarrow\mathbf{(iii)}:$ We have $f^{-1}(0)=g^{-1}(0)$ since  $f\circ\phi=g\circ\phi$ and $\phi^{\ast}$ is injective. Let $\mathfrak{q}=f^{-1}(0)$. For each $s\in S$, by the hypotheses, there is a natural number $n\geq0$ such that $(s/1 +\mathfrak{q}S_{\mathfrak{q}})^{p^{n}}$ is in the image of the induced map $\kappa(\mathfrak{p})\rightarrow\kappa(\mathfrak{q})$ where $\mathfrak{p}=\phi^{\ast}(\mathfrak{q})$ and $p$ is the characteristic of $\kappa(\mathfrak{p})$. Therefore there are elements $r\in R$ and $t\in R\setminus\mathfrak{p}$ such that $\phi(t)s^{p^{n}}-\phi(r)\in\mathfrak{q}$. Thus $f(s)^{p^{n}}=(f\circ\phi)(r)(f\circ\phi)(t)^{-1}=g(s)^{p^{n}}$. This implies that $\big(f(s)-g(s)\big)^{p^{n}}=0$ since the characteristic of $K$ is equal to $p$.  Therefore $f(s)=g(s)$ for all $s\in S$. $\Box$  \\

\end{document}